\documentclass[10pt,twoside]{var_cras_i-fr}
\usepackage[latin1]{inputenc}	
\usepackage[T1]{fontenc}	
\usepackage{amssymb,amsbsy,amsmath,amsfonts,amssymb,amscd}
\usepackage{latexsym,eucal,exscale,epic,eepic,epsfig}
\usepackage[english,francais]{babel}
\usepackage{times}
\input{utile.def}

\renewcommand{\x}{\mathbf{x}}
\renewcommand{\f}{\mathbf{f}}
\renewcommand{\g}{\mathbf{g}}

\renewcommand{\m}{\mathbf{m}}
\newcommand{\gc}{\mathbf{c}}
\newcommand{\V}[1][n]{\cV_{#1}}
\newcommand{\Pol}[1][n]{\KK[\x_{{i,j}}]}  
\newcommand{\Inv}[1][n]{\cI_{#1}} 
\newcommand{\Rec}[1][n]{\cR_{#1}} 
\newcommand{\QInv}[1][n]{\overline\cI_{#1}}
\newcommand{\QRec}[1][n]{\overline\cR_{#1}}
\newcommand{\oast}{{\odot\!\!\!\!\ast}}
\newcommand{\exps}[1]{{\x^{{#1}}}^{\oast}}
\newcommand{\expm}[1]{\x^{#1}}
\newcommand{\upcomplement}{{}^\complement}
\newcommand{\binomial}[2]{\ensuremath{\operatorname{\mathsf C}_{#1}^{#2}}\xspace}
\newtheorem{example}{\it Exemple\rm}[section]

\newtheorem{lemme}[example]{Lemme}
\newtheorem{proposition}[example]{Proposition}
\newtheorem{theorem}[example]{Théorème}

\newtheorem{problem}[example]{Problème}

\newtheorem{conjecture}[example]{Conjecture}

\ISSN{S0764-4442}
\PIT{FLA}
\PXMA{2137}
\Add{1}
\Volume{333}
\Year{2001}
\AuteurCourant{M. Pouzet, N.M. Thiéry}
\TitreCourant{Invariants algébriques de graphes et reconstruction}
\Journal
\Rubrique{Combinatoire}{Combinatorics}
\SousRubrique{}{}
\PresentePar{Gilles}{Kahn}
\Recu{6 septembre 2001}{}{13 septembre 2001}
\begin{document}\label{firstpage}
\selectlanguage{francais}
\title{%
Invariants algébriques de graphes et reconstruction
}
\author{%
Maurice POUZET,\ \ Nicolas M. THIÉRY
}
\address{%
Laboratoire de probabilités, combinatoire et statistique,
Université Claude Bernard Lyon I, \\
Bâtiment recherche [B], 50, avenue Tony-Garnier, Domaine de Gerland, 69366 Lyon Cedex 07\\
Courriel: Maurice.Pouzet@univ-lyon1.fr;
nthiery@users.sf.net
}
\maketitle
\thispagestyle{empty}
\vskip-5mm
\begin{Resume}{%
Nous présentons les résultats d'une étude expérimentale de
l'algèbre des invariants polynomiaux de graphes, motivée par les
problèmes d'isomorphie et de reconstruction.
}\end{Resume}
\selectlanguage{english}
\begin{Etitle}{%
Algebraic graph invariants and reconstruction
}\end{Etitle}
\begin{Abstract}{%
We report on results about a study of algebraic graph invariants,
based on computer exploration, and motivated by graph-isomorphism
and reconstruction problems.
}\end{Abstract}

\AEv

Ulam's reconstruction conjecture \cite{B91} asserts that a graph on
$n$ vertices, $n\geq 3$, is determined up to an isomorphism by the
collection of its vertex-deleted subgraphs. In the frame of invariant
theory, a relevant invariant ring is the ring $\Inv$ of the
polynomials in the indeterminates $x_{{i,j}}$, indexed by pairs
${i,j}$ of ${1,\dots,n}$, which are invariant under the
permutations of ${1,\dots,n}$. A minimal generating set of $\Inv[4]$
composed of $9$ invariants of degree $\leq 5$ was given in
1996 \cite{A_al96}. A computer exploration has been made by the second
author; due to limitations of existing software, a new library
\texttt{PerMuVAR} \cite{T01} was implemented. A minimal generating set
of $\Inv[5]$ composed of $57$ invariants of degree $\leq 9$ was
proposed \cite{T99,T00}, and checked thoroughly later \cite{K00}.
Finally, $\Inv[6]$ seems to be generated by $567$ invariants of degree
$\leq 11$. This exploration led to conjecture that the invariants
$E_1,\dots,E_n$, $p_2,\dots, p_{\binomial{n-1}{2}}$, where
$E_d:=\sum_i\big(\sum_{j\ne i} x_{{i,j}}\big)^d$, and $p_d:=\sum
x_{{i,j}}^d$, form a system of parameters of $\Inv$. This conjecture
was checked up to $n=5$ by Gröbner basis computations; for $n=6$,
the computation was intractable, even on a dedicated parallel
machine \cite{GW99}.

This exploration also leads to several problems and results of
structural nature linked to the reconstruction conjecture on which we
report here. In order to present it, let us fix our notations. Let
$\VV$ be a set. A \emph{weighted graph over $\VV$} is a pair
$\g:=(V(\g),E(\g))$, where $V(\g)$ is the set of \emph{vertices} of
$\g$, and $E(\g)$ is a map from the set of pairs of $V(\g)$ into
$\VV$. If $\VV:={0,1}$ (resp. $V:=\NN$), then $\g$ is a \emph{graph}
(resp. a \emph{multigraph\/}), whose \emph{edges} are pairs ${i,j}$
such that $E(g_{{i,j}})\neq 0$. Denote by $\g_{-i}$ the weighted
graph induced by $\g$ on $V(\g)\smallsetminus{i}$. Two weighted
graphs $\g$ and $\g'$ are \emph{isomorphic} if there is some bijective
map $\sigma:V(\g)\to V(\g')$ such that
$E(g_{{i,j}})=E(g'_{{\sigma(i),\sigma(j)}})$ for all pairs
${i,j}$. They are \emph{hypomorphic} if $\g_{-i}$ and
$\g_{-\sigma(i)}$ are isomorphic for all $i$. Let $n$ be an integer;
$\V(\VV)$ denotes the set of weighted graphs whose vertex set is
${1,\dots,n}$. Let $\WW$ be a set; a map $f:\V(\VV)\to\WW$ is a
\emph{graph invariant\/} (resp. is \emph{reconstructible}) if
$f(\g)=f(\g')$ whenever $\g$ and $\g'$ are isomorphic (resp.
hypomorphic). Let $\KK$ be a field with characteristic $0$; if
$\VV=\WW=\KK$, then the invariant ring $\Inv:=\Pol^{\gS_n}$ identifies
to the subalgebra of polynomial graphs invariants. Given
$\m\in\V(\NN)$, define the monomial $\expm\m:=\prod x_{{i,j}}
^{m_{{i,j}}}$ and the invariant $\exps\m:=\sum_{\m'}\expm{\m'}$,
where $\m'$ runs through the orbit of $\m$. Let $\Rec$ be the
subalgebra of $\Inv$ generated by the $\exps\m$, where $\m$ has at
least one isolated vertex. An invariant $p\in\Rec$ and a multigraph
$\m$ such that $\exps\m\in\Rec$ are called \emph{algebraically
  reconstructible}.

Polynomial invariants separate non-isomorphic weighted graphs.
Similarly, algebraically reconstructible invariants are
reconstructible and separate non-hypomorphic weighted graphs. If $\m$
is a multigraph and $\exps\m$ is reconstructible, then $\m$ is
reconstructible; if moreover $\m$ is a simple graph, then two
hypomorphic simple graphs contain the same number of subgraphs
isomorphic to $\m$. Tutte \cite{T79} proved that $\det(\g)$, the
determinant of the adjacency matrix of $\g$, is reconstructible. The
observation that it is in fact in $\Rec$ led the first
author \cite{P77} to ask whether $\Inv=\Rec$. We show that this is not
the case. Since $\Inv$ and $\Rec$ are graded, it suffices to compare
the dimensions of their homogeneous components of degree $d$; a
computation gives $\dim\Rec[11,18]<\dim\Inv[11,18]$. Therefore, there
exist non-algebraically reconstructible multigraphs with $11$ nodes
and $18$ edges. By the same method, applied to a quotient of $\Inv$,
there exist non-algebraically reconstructible simple graphs with $13$
vertices and $17$ edges, whereas all simple graphs with $13$ vertices
and $17$ edges are reconstructible \cite{MK97}. However, we have no
concrete examples of non-algebraically reconstructible graphs and do
not know if there are some with less vertices and edges.

Among the properties of algebraically reconstructible multigraphs,
note that: (a) if a multigraph $\m$ is not connected then $\m$ is
algebraically reconstructible; (b) if $\m$ is algebraically
reconstructible, then $k\m$ and $k\cdot\gc-\m$ too ($\gc$ is the
complete graph). The reconstructibility, already known or new, of
several graph invariants, readily follows from (a): the maximum size
of a matching, the number of spanning trees, the chromatic number,
the numbers of Hamiltonian paths and cycles, and --- among the new
ones --- the point arboricity, the linear point arboricity, and the
$k$-point partition number.

Characterizing algebraically reconstructible multigraphs seems to
be an interesting problem. We conjecture that \emph{trees are
algebraically reconstructible}. This can be settled in a graded
quotient of $\Inv$, the \emph{forest algebra} $\cF_n$. Contrarily
to the case of $\Inv$, the number of algebraically reconstruc­tible
linear combinations of trees is larger than $\dim\cF_{n,n-1}$, but
we do not know if enough are independent. This leads to problems
about incidence matrices in the vein of Kantor results. Computations
shows that our conjecture holds up to $n=19$. Octopus and stared
octopus are algebraically reconstructible (an \emph{octopus} is a
tree made of paths having a common root; a \emph{stared octopus}
is an octopus with some added stars (trees with diameter $2$)
pending from the root; this includes all trees of diameter at
most $4$). Except for stared multigraphs (all edges share a common
vertex), we do not know if this result extends to multigraphs whose
underlying graph is an octopus. Algebraic versions of the
reconstruction conjecture have been considered by Kocay \cite{K82},
and Mnukhin \cite{M92} ({\it see} also \cite{C96}). The line we
follow here, initiating in \cite{P76,P77}, is parallel to theirs.
The algebra they consider is the quotient $\QInv$ of $\Inv$ by the
ideal generated by the $x_{{i,j}}^2-x_{{i,j}}$. Kocay obtained
the reconstructibility of several parameters by showing that they
are in $\QRec$. He conjectured that trees are in $\QRec$, implying
that hypomorphic simple graphs contains the same number of spanning
trees of a given isomorphism type. He showed that the number of
independent edges-identities provided by members of $\QRec$ can be
smaller than the number of graphs on $n$ vertices and $d$ edges.
However, since $\QInv$ is not graded, this does not imply that
$\QRec\subsetneq\QInv$. In fact, Mnukhin proved that $\QRec=\QInv$
if and only if the reconstruction conjecture holds \cite{M92}.

\par\medskip\centerline{\rule{2cm}{0.2mm}}\medskip
\setcounter{section}{0}
\selectlanguage{francais}

\section{Invariants, graphes et reconstruction, reconstructibilité algébrique}

Soient $\KK$ un corps de caractéristique zéro, $n$ un entier et
$[n]^2$ l'ensemble des paires ${i,j}$ d'éléments de
${1,\dots,n}$. Nous considérons d'une part l'algèbre $\Pol$ des
polynômes en les variables $x_{{i,j}}$ indexées par les
éléments de $[n]^2$ et d'autre part l'ensemble $\V$ des
applications de $[n]^2$ dans $\KK$. Un élément $\g$ de $\V$ est
appelé \emph{graphe valué}. Les \emph{arêtes} de $\g$ sont les
paires ${{i,j}}$ telles que $g_{{i,j}}\neq 0$. Si
$i\in{1,\dots,n}$, nous notons $\g_{-i}$ le graphe obtenu en
supprimant toutes les arêtes de $\g$ incidentes à $i$. Un graphe
qui ne prend que les valeurs $0$ ou $1$ (resp. que des valeurs
entières) s'identifie à un \emph{graphe simple} (resp. un
\emph{multigraphe}). Le groupe symétrique $\gS_n$ agit naturellement
sur $\Pol$ par $\sigma\cdot x_{{i,j}}:= x_{{\sigma(i),\sigma(j)}}$
et sur $\V$ par $\sigma\cdot
g_{{i,j}}:=g_{{\sigma^{-1}(i),\sigma^{-1}(j)}}$. Nous désignons
par $\Inv:=\Pol^{\gS_n}$ la sous-algèbre des polynômes invariants
(c'est-à-dire des polynômes $p$ tels que $\sigma\cdot p=p$ pour
tout $\sigma\in\gS_n$). Deux graphes $\g$ et $\g'$ sont
\emph{isomorphes} s'il existe $\sigma\in\gS_n$ telle que
$\sigma\cdot\g=\g'$. Ils sont \emph{hypomorphes} si $\g_{-i}$ et
$\g_{-\sigma(i)}$ sont isomorphes pour tout $i$. L'\emph{orbite} d'un
graphe valué $\g$ est l'ensemble des graphes qui lui sont
isomorphes. Une telle orbite est appelée \emph{graphe non
  étiqueté}. Sauf ambiguïté, nous laissons le contexte
indiquer si les graphes considérés sont étiquetés ou non. Un
graphe $\g$ est \emph{reconstructible} si tout graphe hypomorphe à
$\g$ lui est isomorphe. Une application $f$ définie sur $\V$ et à
valeurs dans $\KK$ est un \emph{invariant\/} de graphes (valués)
(resp. est \emph{reconstructible}), si elle donne la même valeur à
deux graphes isomorphes (resp. hypomorphes). L'algèbre $\Inv$
s'identifie à l'algèbre des invariants polynomiaux de graphes
valués. À un multigraphe $\m$, nous associons le monôme
$\expm\m:=\prod x_{{i,j}}^{m_{{i,j}}}$ et le polynôme invariant
$\exps\m:=\sum_{\m'}\expm{\m'}$, où $\m'$ parcourt l'orbite de $\m$.
L'ensemble des $\exps\m$ forme une base d'espace vectoriel de $\Inv$.
Nous notons $\Rec$ la sous-algèbre de $\Inv$ engendrée par les
$\exps\m$, où $\m$ parcourt les multigraphes ayant au moins un
sommet isolé. Les éléments de $\Rec$ et les multigraphes $\m$
tels que $\exps\m\in\Rec$ sont dits \emph{algébriquement
  reconstructibles}. Par exemple, la fonction symétrique puissance
$p_k:=\sum x_{{i,j}}^k$ est égale à $\exps{\m}$, où $\m$ est
le multigraphe composé d'une arête valuée $k$; si $n\ge3$, ce
multigraphe a $n-2$ sommets isolés, et donc $p_k\in\Rec$. Il
s'ensuit que tous les polynômes symétriques sont algébriquement
reconstructibles.

Il a été conjecturé par Ulam \cite{B91} que tous les graphes
simples sur au moins $3$ sommets sont reconstructibles. Le lemme
suivant étend à l'hypomorphie un résultat classique, et
suggère une approche de cette conjecture dans le cadre de la
théorie des invariants.

\begin{lemme}\label{lem:separe}
Deux graphes valués dans $\KK$ donnent la même évaluation
à tous les invariants de $\Inv$ {\rm(}resp. de $\Rec${\rm)} si,
et seulement si, ils sont isomorphes {\rm(}resp. hypomorphes{\rm)}.
\end{lemme}

Dans \cite{P77}, le premier auteur demande si $\Inv=\Rec$. La
réponse est négative ({\it cf.\/} § \ref{algrec.limites}).
Toutefois, la considération des polynômes associés aux
multigraphes, et particulièrement aux graphes simples, met en
évidence une propriété forte de reconstructibilité:

\begin{proposition}\label{prop:algrec_rec}
{\rm(a)} Si $\g$ est un graphe simple et $\exps\g$ est
reconstructible, alors $\g$ apparaît le même nombre de fois,
en tant que sous-graphe, dans deux graphes simples hypomorphes;
\par
{\rm(b)} Si $\m$ est un multigraphe et $\exps\m$ est reconstructible,
alors $\m$ est reconstructible.
\end{proposition}

\begin{problem}
Est-ce que {\rm(a)} s'étend aux multigraphes en disant que $\m$
est un \emph{sous-multigraphe} de $\m'$ si pour toute paire
${{i,j}}$ on a $m_{{i,j}}\le m'_{{i,j}}$? Est-ce que la
réciproque de {\rm(a)} est vraie?
\end{problem}

Soit $\m$ un multigraphe. Du lemme \ref{lem:separe} résulte que,
si $\exps\m\in\Rec$, alors $\exps\m$ est reconstructible, et donc,
d'après (b), $\m$ est reconstructible. Si, de plus, $\m$ est un
graphe simple, il satisfait la conclusion de (a). Celle-ci est
apparemment une propriété forte. En effet, les chemins et les
cycles sont trivialement reconstructibles, mais qu'ils apparaissent
le même nombre de fois dans deux graphes simples hypomorphes est
un résultat significatif dû à Tutte \cite{T79}. De même,
les arbres sont reconstructibles \cite{K57}, mais on ne sait pas
si, comme Kocay \cite{K82} l'a conjecturé en 1982, un arbre
apparaît le même nombre de fois dans deux graphes simples
hypomorphes. Les chemins et les cycles sont en fait algébriquement
reconstructibles \cite{P77} ({\it cf.\/} § \ref{algrec.connex}).
Les arbres pourraient être algébriquement reconstructibles
({\it cf.\/} § \ref{algrec.arbres}). En revanche, les graphes
simples ne le sont pas ({\it cf.\/} § \ref{algrec.limites}).

\begin{problem}
Caractériser combinatoirement les graphes, ou les multigraphes,
algébriquement reconstructibles. Quelle est la complexité de ce
problème de reconnaissance? La reconstructibilité algébrique
d'un multigraphe ne dépend-elle que de la partition induite sur
les paires par la valuation?
\end{problem}

\section{Résultats de reconstructibilité découlant de
manipulations algébriques}
\subsection{Non connexité}\label{algrec.connex}. -- 
Soit $\m$ un multigraphe. Soient $\overline\m_1,\dots, \overline\m_k$
les multigraphes sur $n$ sommets obtenus en ajoutant des sommets
isolés aux composantes connexes $\m_1,\dots,\m_k$ de $\m$. Alors,
$\prod^k_{i=1}\exps{\overline\m_i}=\sum^\ell_{i=1}
\alpha_i\exps{\m'_i}$, où les $\alpha_i$ sont des entiers non nuls,
$\m'_1=\m$ et le nombre de composantes connexes des multigraphes
$\m'_i$, $i>1$, est moindre que $k$. De ce simple fait découle:

\begin{proposition}\label{prop:nonconnexes_algrec}
Si $\m$ n'est pas connexe, alors $\exps\m$ est algébriquement
reconstructible.
\end{proposition}

\begin{proposition}\label{prop:bases_ev}
Les polynômes de la forme $\prod^k_{i=1}\exps{\overline\m_i}$ dans
lesquels chaque $\m_i$ est un multigraphe connexe sur $n_i\ge2$
sommets et $\sum^k_{i=1}n_i\leq n$ {\rm(}resp. $n_i\leq n-1${\rm)}
forment une base de l'espace vectoriel $\Inv$ {\rm(}resp. une famille
génératrice de l'espace vectoriel $\Rec${\rm)}.
\end{proposition}

De la proposition \ref{prop:nonconnexes_algrec} découle la
reconstructibilité de plusieurs paramètres attachés aux graphes
simples: la taille maximale d'un couplage, le nombre d'arbres
couvrants, le nombre chromatique \cite{B91}, les paramètres
\emph{point arboricity}, \emph{linear point arboricity} et
\emph{$k$-point partition number} ({\it cf.\/} \cite{H83,LW74} pour
les définitions), le nombre de cycles hamiltoniens, le
déterminant et le polynôme caractéristique \cite{T79}.

\subsection{Dérivation}. -- 
Considérons l'opérateur de dérivation de $\Pol$ défini par
$D:=\sum\partial/\partial x_{{i,j}}$, où la somme est étendue
à toutes les variables $x_{{i,j}}$. Cet opérateur préserve
$\Inv$ et $\Rec$.

\begin{lemme}\label{lem:fractions}
Soit $A$ une sous-algèbre de $\Inv$, préservée par $D$. Si
un polynôme $r$ est quotient de deux polynômes $p$ et $q$ de
$A$ tels que la constante $D^{(\d^{\circ}q)}(q)$ est non nulle,
alors $r$ est élément de $A$.
\end{lemme}

À une constante près non nulle, $D^{(\d^{\circ}q)}(q)$ est la
somme des coefficients de $q$. La conclusion de ce lemme peut être
en défaut lorsque cette somme est nulle. En effet, le corps des
fractions invariantes est engendré par les polynômes
symétriques élémentaires et le polynôme associé au graphe
formé de deux arêtes adjacentes, alors que ces polynômes
n'engendrent pas $\Inv$.

Soit $\m$ un multigraphe et soit $k$ le maximum de $1$ et des
$m_{{i,j}}$. Le \emph{complémentaire de $\m$} est le multigraphe
$\upcomplement\m$ défini par: $\upcomplement m_{{i,j}}
:=k-m_{{i,j}}$. Si $d$ est un entier, on note $d\m$ le multigraphe
$\m'$ défini par $m'_{{i,j}}:=d m_{{i,j}}$. Le
lemme \ref{lem:fractions} a les conséquences suivantes:

\begin{theorem}\label{th:complementaire}
{\rm(a)} $p$ est algébriquement reconstructible si, et seulement
si, $p\prod x_{{i,j}}$ l'est aussi.
\par
{\rm(b)} Si un multigraphe $\m$ est algébriquement reconstructible,
alors $d\m$ l'est aussi.
\par
{\rm(c)} Un multigraphe $\m$ est algébriquement reconstructible
si, et seulement si, $\upcomplement\m$ l'est aussi.
\end{theorem}

\begin{demo}
(a) et (b) ne posent pas de difficultés. (c): soit $\m$ tel que
$\exps\m \in \Rec$, et soit ${\m_1,\dots,\m_l}$ l'orbite de $\m$.
Soit $s_d$ la fonction symétrique puissance de degré $d$ en les
$\expm{\m_i}$. On a $s_d=\exps{d\m}$ et donc, d'après (b),
$s_d\in\Rec$. Ceci étant vrai pour tout $d$, $\Rec$ contient tous
les polynômes symétriques en les $\expm{\m_i}$. Si $e_d$
désigne le polynôme symétrique élémentaire de degré $d$
en les $\expm{\m_i}$, on a $\exps{\upcomplement\m}=
\big(\big(\prod x_{{i,j}}\big)^ke_{\ell-1}\big)/e_\ell$. Donc,
d'après le lemme \ref{lem:fractions}, $\exps{\upcomplement\m}\in\Rec$.
\end{demo}


\begin{problem}
La réciproque de {\rm(b)} est vraie pour les arbres. L'est-elle
pour tout multigraphe?

Pour tout multigraphe $\m$, existe-t-il $d$ tel que $d\m$ soit
algébriquement reconstructible?
\end{problem}

\section{Limites de la reconstructibilité algébrique}
\label{algrec.limites}

Les algèbres $\Inv$ et $\Rec$ sont engendrées par des polynômes
homogènes, et sont donc graduées. La série de Hilbert d'une
algèbre graduée $A$ est $H(A,z):=\sum z^d\dim A_d$, où $A_d$
est la composante homogène de degré $d$ de $A$. Soit $h_{n,d}$
le nombre de multigraphes (non étiquetés) ayant $n$ sommets et
$d$ arêtes, et soit $c_{n,d}$ le nombre de ceux qui sont connexes.
Soit enfin $f_{n,d}$ le nombre de multigraphes à $d$ arêtes,
sans sommets isolés, et dont toutes les composantes connexes ont
au plus $n-1$ sommets. On a $\dim\Inv[n,d]=h_{n,d}$ et $\dim\Rec[n,d]
\le f_{n,d}$. On calcule $H(\Inv, z)$ via une énumération de
Pólya \cite{S79}. Les coefficients $h_{n,d}$ et $f_{n,d}$
s'expriment au moyen de séries génératrices en fonction des
coefficients $c_{n,d}$. Par une méthode analogue à celle
de \cite{HP73}, § 4.2, on en tire les coefficients $f_{n,d}$.

On constate que $\dim\Rec[11,18]\le f_{11,18}<\dim\Inv[11,18]$.
Donc, $\Rec[11]\subsetneq\Inv[11]$, et il existe au moins
$h_{11,18}-f_{11,18}$ multigraphes à $11$ sommets et $18$ arêtes
non algébriquement reconstructibles. Les résultats
expérimentaux pour $n\le19$ suggèrent que, pour tout $n$ et pour
$d$ de l'ordre de $4,1(n-6)$, la proportion de multigraphes
algébriquement reconstructible est majorée par $0,033+6,3
\exp(-0,18 n)$. Un raffinement de ce calcul de dimension dans
l'\emph{algèbre des graphes simples} (quotient de $\Inv$ par
l'idéal engendré par les monômes $x_{{i,j}}^2$) montre qu'il
existe des graphes simples à $13$ sommets et $17$ arêtes non
algébriquement reconstructibles. McKay a vérifié que tous les
graphes simples à $13$ sommets et $17$ arêtes sont
reconstructibles ({\it cf.\/} \cite{MK97}). Il existe donc des
graphes simples reconstructibles et non algébriquement
reconstructibles, mais nous n'en avons pas d'exemple concret.

\section{Reconstructibilité algébrique des arbres}
\label{algrec.arbres}

\begin{conjecture}
Les arbres sont tous algébriquement reconstructibles.
\end{conjecture}

Il suffit de montrer son exactitude dans l'\emph{algèbre des
forêts\/} $\cF_n$ (quotient de l'algèbre des graphes simples
par l'idéal engendré par les cycles). Dans $\cF_n$, les arbres
algébriquement reconstructibles sont engendrés par les forêts
ayant au plus une composante connexe, celle-ci ayant au plus $n-1$
sommets. La dimension de la composante homogène de degré $n-1$
ainsi engendrée est majorée par le nombre total de forêts
non connexes à $n-1$ arêtes. Le quotient de ce nombre par le
nombre d'arbres à $n$ sommets est toujours plus grand que $1$
(il tend vers environ $1,13$). On a donc a priori suffisamment de
générateurs; encore faut-il qu'ils donnent lieu à suffisamment
d'identités linéaires indépendantes. Ceci conduit à un
problème d'algèbre linéaire sur des matrices d'incidence.

Soit $M_n$ (resp. $\overline M_n$) la matrice dont les colonnes sont
indexées par les arbres étiquetés (resp. non étiquetés)
à $n$ sommets et les lignes par les forêts étiquetées (resp.
non étiquetées) à $n-2$ arêtes, le coefficient $a_{\f\g}$
étant égal au nombre d'arêtes de $\g$ dont la suppression
donne $\f$. Pour $M_n$ comme pour $\overline M_n$, le nombre de
lignes est inférieur au nombre de colonnes. On note que
$\overline M_n$ est la transposée de la matrice de
$\exps\f\mapsto\exps\f\sum x_{{i,j}}$ et est équivalente à la
matrice de $\exps\g\mapsto D\exps\g$.

\begin{problem}
{\rm(a)} La matrice $M_n$ est-elle de rang plein?
{\rm(b)} Même question pour $\overline M_n$?
\end{problem}

Une réponse positive pour (a) entraîne une réponse positive
pour (b). Nous avons démontré (a) (resp. (b)) jusqu'à $n=6$
(resp. $n=11$), via le principe d'inclusion-exclusion, et l'avons
contrôlé par le calcul jusqu'à $n=6$ (resp. $n=19$ grâce à
\texttt{linbox} \cite{D01} et \texttt{Nauty} \cite{MK90}, $\overline
M_{19}$ étant de dimension $241029\times317955$). Enfin, \emph{la
  sous-matrice d'incidence des arbres versus les forêts avec un
  sommet isolé est de rang plein}. Les arbres étiquetés étant
les bases du matroïde associé au graphe complet, ce problème
est un cas spécial du suivant: pour quels matroïdes et quels
entiers $p$ et $q$ la matrice d'incidence des parties libres à $q$
éléments versus les parties libres à $p$ éléments est-elle
de rang plein? Quelques résultats \cite{K72,W77,BK92} rentrent dans
ce cadre. Le cas des espaces vectoriels sur les corps finis ne
paraît pas entièrement élucidé.

Une \emph{pieuvre} est un arbre composé de plusieurs chemins issus
d'un même sommet (la tête). Une \emph{pieuvre étoilée} est
une pieuvre à laquelle on a accroché des étoiles (arbres de
diamètre $2$) à la tête. Tous les arbres de diamètre au plus
$4$ sont des pieuvres étoilées.

\begin{lemme}
Soit $p$ un polynôme de l'algèbre engendrée par $\Rec\cup
\Inv[n-1]$. Pour $k\in{1,\dots,n}$, on note $\Theta_k(p):=
\sigma\cdot p$, où $\sigma$ est une permutation quelconque
échangeant $k$ et $n$. Alors, tout polynôme symétrique en
les $\Theta_k(p)$, $k\in{1,\dots,n}$, est algébriquement
reconstructible.
\end{lemme}

\begin{theorem}
Les multigraphes dont toutes les arêtes sont incidentes à un
même sommet, les pieuvres et les pieuvres étoilées sont
algébriquement reconstructibles.
\end{theorem}

\section{Systèmes générateurs de $\Inv$ et $\Rec$}

Une étude de $\Inv$ est faite dans \cite{T99,T00}. Signalons que
pour $n\le4$, le système générateur proposé par Aslaksen
et al. \cite{A_al96} est composé de polynômes invariants
associés à des graphes simples; en outre, $\Inv[4]=\Rec[4]$.
En revanche, pour $n\ge5$, les polynômes invariants associés
aux graphes simples n'engendrent pas l'algèbre des invariants.
Soit $\beta(\Inv)$ le plus petit entier $d$ tel que $\Inv$ soit
engendrée par des polynômes de degré au plus $d$. On a
$\beta(\Inv[3])=3$, $\beta(\Inv[4])=5$, $\beta(\Inv[5])=9$, et
peut-être $\beta(\Inv[6])=11$. Comme la représentation est par
permutation, $\beta(\Inv)\le\binomial{\binomial n2}2$ \cite{GS84}.
La conjecture sur un système de paramètres, mentionnée dans
le résumé, donnerait
$\beta(\Inv)\le\binomial{\binomial{n-1}2}{2}+\binomial{n}{2}$. Enfin,
l'algèbre $\Rec$ est finiment engendrée, et $\beta(\Rec)\leq n
\beta(\Inv[n-1])$.

\Remerciements{Nous remercions J-G. Dumas, P. Flajolet, A. Garsia,
  G. Kemper, B. Kocay, B. McKay et N. Wallach pour leurs nombreuses
  suggestions et pour les résultats de calculs qu'ils nous ont
  communiqués.}
\label{lastpage}
\end{document}